 \documentclass[11pt]{article}

 \usepackage{color}
\usepackage{amsfonts,amssymb,amsmath,mathrsfs}

 \newtheorem{theorem}{Theorem}[section]
 \newtheorem{lemma}[theorem]{Lemma}
 \newtheorem{corol}[theorem]{Corollary}
 \newtheorem{prop}[theorem]{Proposition}
 \newtheorem{example}[theorem]{Example}
\newtheorem{definition}[theorem]{Definition}

 \def\blemma{\begin{lemma}\sl{}\def\elemma{\end{lemma}}}
 \def\bproposition{\begin{prop}\sl{}\def\eproposition{\end{prop}}}
 \def\btheorem{\begin{theorem}\sl{}\def\etheorem{\end{theorem}}}
 \def\bcorollary{\begin{corol}\sl{}\def\ecorollary{\end{corol}}}

 \def\bitemize{\begin{itemize}}\def\eitemize{\end{itemize}}

 \def\beqlb{\begin{eqnarray}}\def\eeqlb{\end{eqnarray}}
 \def\beqnn{\begin{eqnarray*}}\def\eeqnn{\end{eqnarray*}}

 \def\proof{\noindent{\bf Proof.~~}}\def\qed{\hfill$\Box$\medskip}

 \def\<{\langle}\def\>{\rangle}

 \def\lline{---------------------------------}

\begin{document}

\bigskip\bigskip\bigskip
%

\centerline{\LARGE\bf  Asymptotic Behaviors for Critical Branching Processes}
\medskip
\centerline{\LARGE\bf  with  Immigration\footnote{Supported by NSFC
 (No.~11371061)} }

\bigskip

 \centerline{Doudou Li and Mei Zhang\footnote{ Corresponding author}}

 \bigskip

 \centerline{School of Mathematical Sciences, }
 \centerline{Laboratory of Mathematics and Complex Systems,}
 \centerline{ Beijing Normal University,}

 \centerline{Beijing 100875, People's Republic of China}

 \centerline{E-mails: {\tt lidoudou2017@126.com} and {\tt
meizhang@bnu.edu.cn}}
\bigskip

\centerline{\lline\lline\lline}

\medskip

{\narrower

\noindent\textit{\bf Abstract}  In this paper, we investigate the asymptotic behaviors of the critical branching process with immigration $\{Z_n, n\ge 0\}$. First we get some estimation for the probability generating function of $Z_n$. Based on it, we get a large deviation for $Z_{n+1}/Z_n$. Lower and upper deviations for $Z_n$ are also studied. As a by-product, an upper deviation for $\max_{1\le i\le n} Z_i$ is obtained.
\bigskip

\noindent{\it Mathematics Subject Classification (2010):} Primary
60J80; Secondary 60F10

\bigskip
\noindent\textit{\bf Key words and phrases} critical,  branching process, immigration, large deviations.

\medskip

\centerline{\lline\lline\lline}
\begin{math}
~$
\section{Introduction}
\par  Suppose $\{X_{ni},n,i\geq 1\}$
  is a sequence of non-negative integer-valued independent and identically distributed (i.i.d.) random variables with probability generating function $A(x)=\sum_{i=0}^\infty a_ix^i$. \{$Y_n,n\geq 1$\} is another sequence of non-negative integer-valued i.i.d. random variables with probability generating function $B(x)=\sum_{i=0}^\infty b_ix^i$.
$\{X_{ni},~n,i\geq1\}$ are independent with \{$Y_n,n\geq 1$\}. Define $ \{Z_{n}\}$ recursively as
\beqlb\label{1.1}
Z_{n}=\sum_{i=1}^{Z_{n-1}}X_{ni}+Y_n,~~n\geq 1,\quad Z_{0}= 0.
\eeqlb
$\{Z_n, n\ge 0\}$ is called a Galton-Watson branching process with immigration (GMI). Denote $\alpha:=EX_{11}$. When $\alpha>1, \alpha=1~or~\alpha<1$, we shall refer to $\{Z_n\}$ as supercritical, critical and subcritical, respectively. By (\ref{1.1}), the generating function of $Z_n$ can be expressed by
\beqlb\label{1.2}
H_{n}(x)=\prod_{m=0}^{n-1}B[A_{m}(x)],\quad n\ge 1,
\eeqlb
where $A_{m}(x)$ denotes the kth iteration of the function $A(x)$ and $A_{0}(x)=x$.
\par There have been many research works on the large deviations of Galton-Watson branching processes. Particularly, in the critical case, when $Z_{0}=1$ and there is no immigration ($Y_n\equiv 0$), it is known that
\beqlb\label{1.3}
\lim_{n\rightarrow\infty} P\Big\{\bigg |\frac{Z_{n+1}}{Z_{n}}-1\bigg |>\varepsilon \bigg |Z_{n}> 0\Big\}=0.
\eeqlb
Athreya \cite{A99} showed that if $E(Z_{1}^{2r+\delta})<\infty$ for  some $\delta>0$ and $r\geq1$, then  for all $\varepsilon> 0$, there exists $q(\varepsilon)> 0$,
such that
\beqlb\label{1.4}
\lim_{n\rightarrow\infty} nP\Big\{\bigg |\frac{Z_{n+1}}{Z_{n}}-1\bigg |>\varepsilon \bigg |Z_{n}> 0\Big\}=q(\varepsilon)< \infty.
\eeqlb
 In  \cite{N76} and \cite{N06} the authors estimated  the upper   deviation probabilities of $Z_{n}$ and $M_{n}:=\max_{1\le k\leq n}Z_{k}$ under the Cram\'{e}r conditions, respectively. More exactly, in \cite{N76} the inequality
\beqnn
P(Z_{n}\geq k)<  (1+y_{0})\bigg (1+\frac{1}{\frac{1}{y_{0}}+\frac{B_{0}n}{2}} \bigg)^{-k}
\eeqnn
was obtained, where $0< y_{0}< R-1$, $R$ stands for the convergence radius of $A(s)$ and $B_{0}=A^{''}(1+y_{0})$. In \cite{N06} the authors gave that
\beqnn
P(M_{n}\geq k)\leq y_{0}\bigg [\bigg(1+\frac{1}{\frac{1}{y_{0}}+\frac{B_{0}n}{2}}\bigg)^{k} -1\bigg]^{-1}.
\eeqnn

As for the critical Galton-Watson branching processes with immigration, when the functions $A(x)$ and $B(x)$ are analytic in the disk $|x|< 1+\varepsilon$ for some $\varepsilon > 0$, an large deviation  was derived by  \cite{M82}:
\beqnn
\lim_{n\to \infty }P(Z_{n}\geq \frac{bnx}{2})=\frac{1}{\Gamma (\theta)}\int_{x}^{\infty}y^{\theta-1}e^{-y}dy,
\eeqnn
where $b=A^{''}(1-),~\theta =\frac{2B^{'}(1-)}{b}$, $x=o(\frac{n}{\log n})$ and $\Gamma(\cdot)$ is the gamma function.

\par In this paper, we shall study the convergence of the similar type as (\ref{1.4}) for the critical GWI defined by (\ref{1.1}). Some lower deviation probabilities of $Z_{n}$ and upper deviation probabilities of $Z_{n}$ and $M_{n}$ are also established. In the proofs, we have to pay more attention to the changes caused by the immigration and need some precise estimation of the generating function of $Z_n$.

 We will begin our discussion under the following assumption:
\par (H)~$0< a_{0},~b_{0}< 1,~\sum\limits_{j=1}^{\infty}a_{j}j^{2}\log j< \infty,~\sum\limits_{j=1}^{\infty}b_{j}j^{2}<  \infty,~\alpha=1, ~0<\beta:=B^{'}(1-)<\infty,~0< \gamma=\frac{1}{2}A^{''}(1-)< \infty$.

 In the following, we define  $\sigma =\frac{\beta}{\gamma}$. We write   $d_{n}=O(e_n)$ if and only if  there exist $C_{1}$ and $C_{2}$ such that
$$
C_{1}\leq \varliminf\limits_{n\rightarrow \infty}\frac{d_{n}}{e_n}\leq \varlimsup\limits_{n\rightarrow \infty}\frac{d_{n}}{e_n}\leq C_{2};
$$
$d_n\sim e_n$ if and only if
$$
\lim_{n\rightarrow \infty}\frac{d_{n}}{e_n}=1.
$$
$C_1, C_2,\cdots$ are positive constants whose value may vary from place to place.

The rest of the paper is organized as follows. Some preliminary results are given in Section 2.
In Section 3 we state the main theorems. Section 4 is devoted to the proofs of the main theorems.

\section{Preliminary results}

\blemma\label{l2.1}(Athreya and Ney \cite{A72})
 Assume $\alpha=1$, $0< \gamma < \infty$ and let $\delta (x)=\gamma-[\frac{1}{1-A(x)}-\frac{1}{1-x}]$. Define $h_{n}(x)=\sum\limits_{m=0}^{n-1}\delta(A_{m}(x))$ for $n\geq 1$  and $h_{0}(x)= 0$. Then
\beqlb\label{2.1}
\frac{1}{1-x}+n\gamma -\frac{1}{1-A_{n}(x)}=h_{n}(x),~~~0\leq x< 1.
\eeqlb
Furthermore, $\delta(x)$ satisfies the inequality
\beqnn
\frac{-\gamma^{2}(1-x)}{1-a_{0}}\leq \delta (x)\leq \varepsilon (x),~~~0\leq x< 1,
\eeqnn
where $0\leq \varepsilon (x):=\gamma - \frac{A(x)-x}{(1-x)^{2}}$, which is non-increasing in $x$ and $\varepsilon (x)\downarrow 0$ as $x \uparrow 1$.
\elemma

\blemma\label{l2.2} (Pakes \cite[Theorems 1,2]{P72} )  Under condition (H), we have \beqlb\label{4.13}
\lim\limits_{n\rightarrow \infty}n^{\sigma}H_{n}(x)=U(x),
\eeqlb where $U(x)$ satisfies the functional equation
$$
B(x)U(A(x))=U(x).
$$
The above convergence is uniform over compact subsets of the open unit disc. Moreover, \beqlb\label{4.15}
 U(x)\sim (1-x)^{-\sigma},\quad x\to 1^-.
\eeqlb Denoting the power series representation of $U(x)$ by $\sum_{j=0}^\infty \mu_jx^j$, then
\beqlb\label{4.14}
\lim_{n\rightarrow\infty} n^{\sigma}P\{Z_{n}=j\}=\mu_{j}, \quad j\ge 0.
\eeqlb
\elemma

\blemma\label{l2.3} (Pakes \cite[Theorem 10]{P75} )
Let  $p_{0j}^{(n)}$ be the n-step transition probability of $\{Z_n\}$ from state $0$ to $j$ and $~\nu_{n}=\sum\limits_{j=1}^{\infty}\frac{p_{0j}^{(n)}}{j}$.
 Under condition (H),

\noindent (i)~if $\sigma< 1$, then
\beqnn
\nu_{n}\sim n^{-\sigma}\int_{0}^{1}\frac{U(s)-U(0)}{s}ds,
\eeqnn
where $U(s)$ is defined by (\ref{4.13}).

\noindent (ii)~if $\sigma> 1$, then
\beqnn
\nu_{n}\sim \frac{1}{n(\beta-\gamma)}.
\eeqnn

\elemma

\section{Main results}

\btheorem\label{t3.1}
 Assume (H) holds.  For each $\varepsilon >0$, define
\beqlb\label{3.7}
A(k,\varepsilon)= P(|\bar{X}_k+\frac{Y_1}{k}-1|>\varepsilon),
\eeqlb
where $\bar{X}_k=\frac{1}{k}\sum_{i=1}^k X_{1i}$. For $r>\sigma$, if there exists $C_{\varepsilon}> 0$ such that $A(k,\varepsilon)\leq C_\varepsilon  k^{-r}\;$ for all $k\geq1$, then  there exists $q(\varepsilon)> 0$, such that
 \beqlb\label{3.8}
\lim_{n\rightarrow\infty} n^{\sigma}P\Big\{\bigg |\frac{Z_{n+1}}{Z_{n}}-1\bigg |>\varepsilon \bigg |Z_{n}> 0\Big\}=q(\varepsilon)< \infty.
\eeqlb
 \etheorem
\bcorollary
Assume (H) holds, $E(X_{11}^{2r+\delta})<\infty$ and $E(Y_{1}^{r})<\infty$ for  some $\delta>0$ and $r> \max\{\sigma, 1\}$. Then (\ref{3.8}) holds.
\ecorollary

\btheorem\label{t3.3}
Define $J_{n}=D\{\frac{Z_{n+1}}{Z_{n}}|Z_{n}>0\}$.  Assume $0< DY_{1}< \infty$ and  (H) holds. We have

\noindent(1)~if $\sigma < 1$, then
\beqnn
J_{n}=\kappa n^{-\sigma}(1+o(1)),
\eeqnn
where $\kappa =2\gamma \int_{0}^{1}\frac{U(s)-U(0)}{s}ds+D(Y_{1})\sum\limits_{k\geq 1}\frac{\mu_{k}}{k^{2}}$ with $\{\mu_{k}\}$ given by (\ref{4.14});

\noindent (2)~if $\sigma = 1$, then
\beqlb\label{3.9}
J_{n}=O(\frac{\log n}{n});
\eeqlb
(3)~if $\sigma > 1$ and $\sigma \neq 2$, then
\beqnn
J_{n}=\frac{2\gamma}{n(\beta-\gamma)}(1+o(1)).
\eeqnn

\etheorem

\btheorem\label{t3.4}
Assume  (H) holds. Let $k_{n}\to \infty$ and $k_{n}=o(n)$ as $n\to \infty$. Then
\beqnn
P(Z_{n}\leq k_{n})\leq C_{3}(1+\gamma \frac{n}{k_{n}})^{-\sigma},
\eeqnn
 as $n\rightarrow \infty$.
\etheorem
\btheorem\label{t3.5}
Assume (H) holds. Let $R$ stand for the convergence radius of  $A(x)$. Assume $R> 1$, $\frac{k_n}{n}\rightarrow \infty$ and $k_n=o(n^{2})$ as $n\rightarrow \infty$. Then
 \beqlb\label{3.10}
 P(Z_{n}\geq k_n)\leq (\frac{k_n}{\gamma n})^{B'(1+\frac{k_n}{\gamma^{2}n^{2}}-\frac{1}{\gamma n})\frac{1}{\gamma}}\exp\Big \{-\frac{k_n}{\gamma n}+1-\frac{1}{\gamma}\lambda \frac{k_n}{n^{2}}\ln \frac{k_n}{n}\Big \}\bigg (1+O(\frac{k_n}{n^{2}})\bigg )
 \eeqlb
as $n\rightarrow \infty$, where $\lambda =1-\frac{\rho}{6\gamma^{2}}$ and $\rho =A^{'''}(1-)< \infty$.
\etheorem

\bcorollary
Assume the hypotheses of Theorem~\ref{t3.5} hold. Let $M_{n}=\max_{1\le k\leq n}Z_{k}$. Then
\beqlb\label{3.11}
 P(M_{n}\geq k_n)\leq (\frac{k_n}{\gamma n})^{B'(1+\frac{k_n}{\gamma^{2}n^{2}}-\frac{1}{\gamma n})\frac{1}{\gamma}}\exp\Big \{-\frac{k_n}{\gamma n}+1-\frac{1}{\gamma}\lambda \frac{k_n}{n^{2}}\ln \frac{k_n}{n}\Big \}\bigg (1+O(\frac{k_n}{n^{2}})\bigg ).
 \eeqlb
\ecorollary
\noindent{\bf Remark:}  The right sides of (\ref{3.10}) and (\ref{3.11}) approximate to $(\frac{k_n}{\gamma n})^{\sigma}\exp \{-\frac{k_n}{\gamma n}\}$ as $n\to \infty$.

\section{Proofs of main results}
In this section we prove Theorem 3.1-Corollary 3.6. First we present the following proposition.
\bproposition\label{p4.1}
Assume  (H) holds. Then for each $C_4>0$, there exist positive constants $C_{5}$ and $C_{6}$ such that for any   $0< s\leq C_{4}n$,
\beqnn
C_5(1+\gamma s)^{-\sigma}\le H_{n}(e^{-\frac{s}{n}})\le  C_6(1+\gamma s)^{-\sigma}.
\eeqnn
\eproposition

{\proof }
By Taylor's formula, we know that
\beqnn
\log x=x-1-\frac{1}{2\theta_{1}^{2}}(x-1)^{2},~~x\leq \theta_{1}\leq 1, x\in (0,1),
\eeqnn
and
\beqnn
1-B(x)=\beta (1-x)-\frac{B^{''}(\theta_{2})}{2}(1-x)^{2},~~x\leq \theta_{2}\leq 1, x\in (0,1).
\eeqnn
Recalling (\ref{1.2}), we obtain
\beqnn
\log H_{n}(x)
&=&\sum\limits_{m=0}^{n-1}[B(A_{m}(x))-1-\frac{1}{2\theta_{3}^{2}}(B(A_{m}(x))-1)^{2}]\\
&=&-\sum\limits_{m=0}^{n-1}[1-B(A_{m}(x))]-\frac{1}{2\theta_{3}^{2}}\sum\limits_{m=0}^{n-1}(B(A_{m}(x))-1)^{2}\\
&=&-\sum\limits_{m=0}^{n-1}[\beta (1-A_{m}(x))-\frac{B^{''}(\theta_{4})}{2}(1-A_{m}(x))^{2}]-\frac{1}{2\theta_{3}^{2}}\sum\limits_{m=0}^{n-1}(B(A_{m}(x))-1)^{2}\\
&=&-\beta \sum\limits_{m=0}^{n-1}(1-A_{m}(x))+I_{0}(n),
\eeqnn
where \beqnn & & I_{0}(n)=\sum\limits_{m=0}^{n-1}\frac{B^{''}(\theta_{4})}{2}(1-A_{m}(x))^{2}
-\frac{1}{2\theta_{3}^{2}}\sum\limits_{m=0}^{n-1}(B(A_{m}(x))-1)^{2},\\
& & B(A_{m}(x))\leq \theta_{3}\leq 1,~A_{m}(x)\leq \theta_{4}\leq 1,~x\in (0,1).\eeqnn  Since $|1-A_{m}(x)|\leq 2|1-A_{m}(0)|\sim \frac{2}{m\gamma}$ as $m\rightarrow \infty$ (see \cite{P72} Page 74), it is easy to show that $I_{0}(n)$ is uniformly bounded for all $x\in (0,1)$ as $n\rightarrow \infty$.
\par It is known from (\ref{2.1}),
\beqnn
1-A_{m}(x)=\frac{1-x}{1+\gamma m(1-x)}+\frac{1-x}{1+\gamma m(1-x)}\Big [\frac{h_{m}(x)}{\frac{1+\gamma m(1-x)}{1-x}-h_{m}(x)}\Big ].
\eeqnn
Consequently,
\beqnn
& & \log H_{n}(x)\\
&=&-\beta \sum\limits_{m=0}^{n-1}(1-A_{m}(x))+I_{0}(n)\\
&=&-\beta \sum\limits_{m=0}^{n-1}\frac{1-x}{1+\gamma m(1-x)}
-\beta \sum\limits_{m=0}^{n-1}\frac{1-x}{1+\gamma m(1-x)}\Big [\frac{h_{m}(x)}{\frac{1+\gamma m(1-x)}{1-x}-h_{m}(x)}\Big ]+I_{0}(n)\\
&=&-\beta \sum\limits_{m=0}^{n-1}\frac{1-x}{1+\gamma m(1-x)}+I_{1}(n)+I_{0}(n),
\eeqnn
where $I_{1}(n)=-\beta \sum\limits_{m=0}^{n-1}\frac{1-x}{1+\gamma m(1-x)}\Big [\frac{h_{m}(x)}{\frac{1+\gamma m(1-x)}{1-x}-h_{m}(x)}\Big ]$.  By  \cite[Theorem 1]{P72},
\beqnn
\sum\limits_{m=1}^{\infty}\frac{|h_{m}(x)|}{m^{2}}\leq \max \Big\{\frac{\gamma^{2}}{1-A(0)}\sum\limits_{m=1}^{\infty}\frac{1}{m^{2}}\sum\limits_{k=0}^{m-1}(1-A_{k}(0)),\sum\limits_{m=1}^{\infty}\frac{1}{m^{2}}\sum\limits_{k=0}^{m-1}\varepsilon (A_{k}(0))\Big\}< \infty.
\eeqnn
Hence, $\sum\limits_{m=1}^{\infty}\frac{|h_{m}(x)|}{m^{2}}$ is uniformly bounded for all $x\in (0,1)$. Furthermore, we have
\beqnn
& & \Big |\sum\limits_{m=1}^{n-1}\frac{1-x}{1+\gamma m(1-x)}\Big [\frac{h_{m}(x)}{\frac{1+\gamma m(1-x)}{1-x}-h_{m}(x)}\Big ]\Big | \\
&&\leq \sum\limits_{m=1}^{n-1}\frac{1-x}{1+\gamma m(1-x)}\Big |\frac{\frac{h_{m}(x)}{m}}{\frac{1}{m(1-x)}+\gamma-\frac{h_{m}(x)}{m}}\Big |\\
  &&\leq \sum\limits_{m=1}^{n-1}\frac{1}{\gamma m}\Big |\frac{\frac{h_{m}(x)}{m}}{\frac{1}{m(1-x)}+\gamma-\frac{h_{m}(x)}{m}}\Big |\\
 &&\leq \frac{1}{\gamma a_{0}}\sum\limits_{m=1}^{\infty}\frac{|h_{m}(x)|}{m^{2}}.
\eeqnn
Then $I_{1}(n)$ is uniformly bounded for all $x\in (0,1)$ as $n\rightarrow \infty$. Finally, we have
\beqnn
\log H_{n}(x)=-\beta \sum\limits_{m=0}^{n-1}\frac{1-x}{1+\gamma m(1-x)}+O(1)
\eeqnn
 uniformly for $x\in (0,1)$. Let $0<s\le C_4n$ and $x=e^{-\frac{s}{n}}$. Then
\beqlb\label{2.5}
\log H_{n}(e^{-\frac{s}{n}})=-\beta \sum\limits_{m=0}^{n-1}\frac{1-e^{-\frac{s}{n}}}{1+\gamma m(1-e^{-\frac{s}{n}})}+O(1).
\eeqlb
It can be easily observed that
\beqnn
0\leq \frac{\frac{s}{n}}{1+\gamma m \frac{s}{n}}-\frac{1-e^{-\frac{s}{n}}}{1+\gamma m(1-e^{-\frac{s}{n}})}\leq \frac{\frac{s^{2}}{2n^{2}}}{1+\gamma^{2} m^{2} \frac{s}{n}(1-e^{-\frac{s}{n}})}.
\eeqnn
Now we prove there exists $C_{7}$, such that
\beqlb\label{2.6}
0\leq I_{2}(n,s):=\sum\limits_{m=0}^{n-1}\frac{\frac{s^{2}}{2n^{2}}}{1+\gamma^{2} m^{2} \frac{s}{n}(1-e^{-\frac{s}{n}})}\leq C_{7},\quad n\rightarrow\infty.
\eeqlb
To see this, setting $u(t)=\frac{t^{2}}{1+\gamma^{2} m^{2}t(1-e^{-t})}$. Then $u(t)$ is increasing for $t>0$. Hence by (\ref{2.6}) we have
\beqnn
 I_{2}(n,s)\leq \frac{1}{2}\sum\limits_{m=0}^{\infty}\frac{C_4^2}{1+C_4\gamma^{2} m^{2}(1-e^{-C_4})}:=C_7< \infty.
\eeqnn
 Recalling (\ref{2.5}) we obtain
\beqnn
\log H_{n}(e^{-\frac{s}{n}})=-\beta \sum\limits_{m=0}^{n-1}\frac{\frac{s}{n}}{1+\gamma m \frac{s}{n}}+O(1).
\eeqnn
Since
\beqnn
\int_{0}^{s}\frac{1}{1+\gamma x}dx \leq \sum\limits_{m=0}^{n-1}\frac{\frac{s}{n}}{1+\gamma m \frac{s}{n}}\leq \int_{0}^{s}\frac{1}{1+\gamma x}dx+\frac{s}{n}-
\frac{s}{n(1+\gamma s)},
\eeqnn
 we arrive at
\beqnn
-\sigma \log (1+\gamma s)+O(1)\leq \log H_{n}(e^{-\frac{s}{n}})\leq -\sigma \log (1+\gamma s)+O(1).
\eeqnn
The proof is now complete.~\qed

\noindent{\bf Proof of Theorem 3.1.}
Using the branching property, we have
\beqlb\label{4.12}
n^{\sigma}P\Big\{\bigg |\frac{Z_{n+1}}{Z_{n}}-1\bigg |>\varepsilon \bigg |Z_{n}> 0\Big\}
=\sum\limits_{j=1}^{\infty}A(j,\varepsilon)n^{\sigma}P\{Z_{n}=j|Z_{n}> 0\},
\eeqlb
where $A(j,\varepsilon)$ is given by (\ref{3.7}).
\par From (\ref{4.14}), $\lim\limits_{n\rightarrow\infty} P\{Z_{n}=0\}=0$, then the condition on $Z_{n}> 0$ is not necessary when we consider the case $n\rightarrow\infty$. Therefore, in the following we only consider
\beqnn
n^{\sigma}P\Big\{\bigg |\frac{Z_{n+1}}{Z_{n}}-1\bigg |>\varepsilon \Big\}
=\sum\limits_{j=1}^{\infty}A(j,\varepsilon)n^{\sigma}P\{Z_{n}=j\}.
\eeqnn
Next, we will prove as $n\rightarrow\infty$,
\beqlb\label{4.16}
 n^{\sigma}P\Big\{\bigg |\frac{Z_{n+1}}{Z_{n}}-1\bigg |>\varepsilon \Big\}\rightarrow \sum\limits_{j=1}^{\infty}A(j,\varepsilon)\mu_{j}<\infty.
 \eeqlb
Since $A(j,\varepsilon)\leq C_{\varepsilon}j^{-r}$, and $j^{-r}\sim (j+1)^{-r}$ as $j\rightarrow \infty$, then there exists $C_{\varepsilon}'$ such that $A(j,\varepsilon)\leq C_{\varepsilon}'(j+1)^{-r}$ for all $j\ge 1$. Therefore,
\beqnn
l_{n}(j):=n^{\sigma}A(j,\varepsilon)P\{Z_{n}=j\}\leq n^{\sigma}C_{\varepsilon}'(j+1)^{-r}P\{Z_{n}=j\}:=\tilde{l}_{n}(j).
\eeqnn
Using (\ref{4.14}), we have for $j\ge 0$,
\beqnn
\lim\limits_{n\rightarrow \infty}\tilde{l}_{n}(j)=C_{\varepsilon}'\mu_{j}(j+1)^{-r}:=\tilde{l}(j).
\eeqnn
 By \cite{M82b},
\beqlb\label{4.22}
  \mu_{j}\sim (\gamma^{\sigma}\Gamma(\sigma))^{-1}j^{\sigma-1},\quad j\rightarrow \infty.
 \eeqlb
Then for $r> \sigma$,
 $$\sum\limits_{j=0}^{\infty}\tilde{l}(j)=C_{\varepsilon}'\sum\limits_{j=0}^{\infty}\mu_{j}(j+1)^{-r} < \infty.$$
Now, using a  modification of the Lebesgue dominated convergence theorem, it is sufficient to show that as $n\to \infty$,
\beqlb\label{lnj}
\sum\limits_{j=0}^{\infty}\tilde{l}_{n}(j)\to \sum\limits_{j=0}^{\infty}\tilde{l}(j),
\eeqlb
which is equivalent to
\beqlb\label{4.19a}
\lim\limits_{n\rightarrow \infty} n^{\sigma}E((Z_{n}+1)^{-r})=\sum\limits_{j=0}^{\infty}\mu_{j}(j+1)^{-r}.
\eeqlb
In the following we prove (\ref{4.19a}). For $r>0$, we have
\beqlb\label{harmo}
\Gamma (r)n^{\sigma}E((Z_{n}+1)^{-r})
&=& \int_{0}^{\infty}n^{\sigma}E(e^{-t(Z_{n}+1)})t^{r-1}dt\nonumber\\
&=& \int_{0}^{1}n^{\sigma}H_{n}(s)(-\log s)^{r-1}ds\nonumber\\
&=& I_{3}(n)+I_{4}(n)
\eeqlb
where $$I_{3}(n)=\int_{0}^{\frac{1}{e}}n^{\sigma}H_{n}(s)(-\log s)^{r-1}ds,$$
and $$I_{4}(n)=\int_{\frac{1}{e}}^{1^{-}}n^{\sigma}H_{n}(s)(-\log s)^{r-1}ds. $$
It is easy to see $\int_{0}^{\frac{1}{e}}(-\log s)^{r-1}ds < \infty$.  By Lemma~\ref{l2.2}, $U(s)$ is bounded in $[0, \frac{1}{e}]$. Therefore
\beqlb\label{4.17}
\lim\limits_{n\rightarrow \infty}I_{3}(n)=\int_{0}^{\frac{1}{e}}U(s)(-\log s)^{r-1}ds< \infty.
\eeqlb
Define $$f_{n}(s)=n^{\sigma}H_{n}(s)(-\log s)^{r-1},\quad  s\in (0,1), $$ then $$f_{n}(s)\rightarrow f(s):=U(s)(-\log s)^{r-1}. $$
From Proposition~\ref{p4.1}, we know that for $t\in [e^{-1},1)$, there exists $N$ and $C_{8}$ such that for $n>N$,
\beqnn
H_{n}(t)\leq C_{8}(1-\gamma n\log t)^{-\sigma},
\eeqnn
hence,
\beqnn
f_{n}(s)\leq C_{8}n^{\sigma}(1-\gamma n\log s)^{-\sigma}(-\log s)^{r-1}:=g_{n}(s)
\eeqnn
for all $n>N$. It is not difficult to see
\beqnn
g_{n}(s)\nearrow C_{8}(-\gamma \log s)^{-\sigma}(-\log s)^{r-1}:=g(s),
\eeqnn
and for $r> \sigma$,
\beqnn
\int_{\frac{1}{e}}^{1}g(s)ds=C_{8}\int_{0}^{1}(\gamma t)^{-\sigma}t^{r-1}e^{-t}dt< \infty.
\eeqnn
Using the  modification of  dominated convergence theorem, we have
\beqlb\label{4.19}
\int_{\frac{1}{e}}^{1}f_{n}(s)ds\longrightarrow \int_{\frac{1}{e}}^{1}f(s)ds,\quad n\rightarrow \infty.
\eeqlb
By a change of variable $u=-\log s$, the right side of (\ref{4.19}) turns out to be
\beqnn
\int_{0}^{1}U(e^{-u})u^{r-1}e^{-u}du,
\eeqnn
 which is finite by using (\ref{4.15}). Hence, we obtain
\beqlb\label{4.20}
\lim\limits_{n\rightarrow \infty}I_{4}(n)=\int_{\frac{1}{e}}^{1}f(s)ds< \infty.
\eeqlb
Together with (\ref{harmo})--(\ref{4.17}) and (\ref{4.20}), we have
\beqnn
\Gamma (r)n^{\sigma}E((Z_{n}+1)^{-r})\to\int_{0}^{1}f(s)ds<\infty, \quad n\to \infty,
\eeqnn
which yields
\beqnn
\lim\limits_{n\rightarrow \infty}\sum\limits_{j=0}^{\infty}\tilde{l}_{n}(j)=\frac{C'_{\varepsilon}}{\Gamma (r)}\int_{0}^{1}f(s)ds<\infty.
\eeqnn
Clearly, $$
\frac{C'_{\varepsilon}}{\Gamma (r)}\int_{0}^{1}f(s)ds=\sum\limits_{j=0}^{\infty}\tilde{l}(j).$$
Thus we get (\ref{lnj}), and then(\ref{4.16}) holds. The proof is completed.~\qed

 \noindent{\bf Proof of Corollary 3.2.}
 By Markov's inequality, we have
 $$A(k,\varepsilon)= P(|\bar{X}_k+\frac{Y_1}{k}-1|>\varepsilon)\leq \frac{E(\sqrt{k}(\bar{X}_k+\frac{Y_1}{k}-1))^{2r}}{\varepsilon^{2r}k^{r}}.$$
 Using the assumption and \cite[Page 112, section 9.9]{L10}, we obtain
 $$\widetilde{C_{\varepsilon}}=\sup\limits_{k}E(\sqrt{k}(\bar{X}_k+\frac{Y_1}{k}-1))^{2r}< \infty.$$
 Then there exists a constant $C_{\varepsilon}$ such that $A(k,\varepsilon)\leq C_{\varepsilon}k^{-r}$ for all $k\geq 1$.
 ~\qed

\noindent{\bf Proof of Theorem 3.3.}
 Let $p_{nk}^{\ast}=P(Z_{n}=k|Z_{n}> 0)$. By direct calculation, we have
\beqlb\label{4.21a}
D\{\frac{Z_{n+1}}{Z_{n}}|Z_{n}> 0\}=D(X_{11})\sum\limits_{k\geq 1}\frac{p_{nk}^{\ast}}{k}+D(Y_{1})\sum\limits_{k\geq 1}\frac{p_{nk}^{\ast}}{k^{2}}.
\eeqlb
First we discuss $\sum\limits_{k\geq 1}\frac{p_{nk}^{\ast}}{k}$. The cases $\sigma\neq1$ have been given by Lemma~\ref{l2.3}. For $\sigma=1$, we shall prove \beqlb\label{11}\sum\limits_{k\geq 1}\frac{p_{nk}^{\ast}}{k}=O(\frac{\log n}{n}). \eeqlb We know that
\beqnn
 \sum\limits_{k\geq 1}\frac{p_{nk}^{\ast}}{k}=\int_{0}^{1}\frac{H_{n}(x)-H_{n}(0)}{x(1-H_{n}(0))}dx:=I_{5}(n)+I_{6}(n),
 \eeqnn
 where $$I_{5}(n)=\int_{0}^{e^{-\frac{s}{n}}}\frac{H_{n}(x)-H_{n}(0)}{x(1-H_{n}(0))}dx, $$
 and  $$I_{6}(n)=\int_{e^{-\frac{s}{n}}}^{1}\frac{H_{n}(x)-H_{n}(0)}{x(1-H_{n}(0))}dx, $$ with $s=\frac{n}{\log n}$. Using Lemma~\ref{l2.2} we have that
 \beqnn
 \int_{e^{-\frac{s}{n}}}^{1}\frac{H_{n}(0)}{x(1-H_{n}(0))}dx=\frac{H_{n}(0)}{\log n(1-H_{n}(0))}\sim \frac{\mu_{0}}{n^{\sigma}\log n},\quad n\rightarrow \infty.
 \eeqnn
 Let $I_{6}^{'}(n)=\int_{e^{-\frac{s}{n}}}^{1}\frac{H_{n}(x)}{x}dx$. Next we consider the order of $I_{6}^{'}(n)$. By Proposition~\ref{p4.1},
 \beqnn
 I_{6}^{'}(n)=\frac{1}{n}\int_{0}^{\frac{n}{\log n}}H_{n}(e^{-\frac{\theta}{n}})d\theta =O(\frac{\log n}{n}).
 \eeqnn
Moreover, noticing that
$$
\frac{H_{n}(x)-H_{n}(0)}{ x(1-H_{n}(0))}=E(x^{Z_n-1}|Z_n>0)
$$
is non-decreasing in $x$, then by the definition of $I_{5}(n)$ and Proposition 4.1, we obtain  \beqnn
 I_{5}(n)\leq \frac{H_{n}(e^{-\frac{s}{n}})}{1-H_{n}(0)}=O(\frac{\log n}{n}).
 \eeqnn
Thus (\ref{11}) holds.

Now we turn to estimate $$\nu_{n}^{\ast}:=\sum\limits_{k\geq 1}\frac{p_{nk}^{\ast}}{k^{2}}. $$
(i)~if $\sigma < 2$, by (\ref{4.16}), (\ref{4.14}) and (\ref{4.22}),
 \beqnn
 \lim\limits_{n\rightarrow \infty}n^{\sigma}\nu_{n}^{\ast}=\sum\limits_{k=1}^{\infty}\lim\limits_{n\rightarrow \infty}\frac{n^{\sigma}p_{nk}^{\ast}}{k^{2}}=
 \sum\limits_{k=1}^{\infty}\frac{\mu_{k}}{k^{2}}<\infty.
 \eeqnn
 Then we have
  \beqlb\label{4.23a}
  \nu_{n}^{\ast}\sim n^{-\sigma}\sum\limits_{k\geq 1}\frac{\mu_{k}}{k^{2}},\quad n\rightarrow \infty.
  \eeqlb
 (ii)~if $\sigma > 2$, first it is known that
 \beqnn
 \Gamma (2)n^{2}\nu_{n}^{\ast}
 &=& \int_{0}^{\infty}n^{2}E(e^{-tZ_{n}}|Z_{n}>0)tdt\\
 &=& \frac{1}{P(Z_{n}>0)}\bigg(\int_{0}^{1}n^{2}E(e^{-tZ_{n}}, Z_{n}>0)tdt+\int_{1}^{\infty}n^{2}E(e^{-tZ_{n}}, Z_{n}>0)tdt\bigg)\\
 &:=& \frac{1}{P(Z_{n}>0)}(I_{7}(n)+I_{8}(n)).
 \eeqnn
  By a change of variable $t=\frac{s}{n}$, we have
\beqnn
I_{7}(n)
&=& \int_{0}^{n} (H_{n}(e^{-\frac{s}{n}})-H_{n}(0))sds\\
&=& \int_{0}^{\infty}I_{(s\leq n)} (H_{n}(e^{-\frac{s}{n}})-H_{n}(0))sds\\
&:=& \int_{0}^{\infty}q_{n}(s)sds.
\eeqnn
Using Proposition~\ref{p4.1}, there exists $C_{9}$ such that
\beqnn
q_{n}(s)\leq C_{9}(1+\gamma s)^{-\sigma} :=v(s).
\eeqnn
Clearly, for  $\sigma >2$,
\beqnn
\int_{0}^{\infty}v(s)sds< \infty.
\eeqnn
 Therefore, using the dominated convergence theorem, we have
\beqnn
\lim\limits_{n\rightarrow \infty}I_{7}(n)=\int_{0}^{\infty}q(s)sds < \infty,
\eeqnn
where  by \cite[Theorem 3]{P69},
$$q(s):=\lim\limits_{n\rightarrow \infty}q_{n}(s)=(1+\gamma s)^{-\sigma}. $$

For $I_8(n)$,
\beqnn
I_{8}(n)
&=& \int_{1}^{\infty}n^{2}E(e^{-tZ_{n}} I_{(Z_{n}=1)})tdt +\int_{1}^{\infty}n^{2}E(e^{-tZ_{n}}I_{(Z_{n}\geq 2)})tdt\\
&:=& K_{1}(n)+K_{2}(n).
\eeqnn
By (\ref{4.14}) and the Lebesgue dominated convergence theorem,
\beqnn
K_{1}(n)=\int_{1}^{\infty}n^{2}P(Z_{n}=1)e^{-t}tdt\rightarrow 0,\quad n\rightarrow \infty.
\eeqnn
Meanwhile, by (\ref{4.13}) and the Lebesgue dominated convergence theorem,
\beqnn
K_{2}(n)
&=& \int_{1}^{\infty}n^{2}E(e^{-t(Z_{n}-1)}I_{(Z_{n}\geq 2)})e^{-t}tdt\\
&\leq & \int_{1}^{\infty}n^{2}E(e^{-\frac{t}{2}Z_{n}})e^{-t}tdt\\
&=& \int_{1}^{\infty}n^{2}H_{n}(e^{-\frac{t}{2}})e^{-t}tdt\rightarrow 0,\quad n\rightarrow \infty.
\eeqnn
Then, we get
\beqnn
\lim\limits_{n\rightarrow \infty}I_{8}(n)=0.
\eeqnn
Let $C_{10}=\int_{0}^{\infty}q(s)sds$, we get
 \beqlb\label{4.24}
 \nu_{n}^{\ast}\sim C_{10}n^{-2},\quad n\rightarrow \infty.
 \eeqlb
 Collecting (\ref{4.21a})--(\ref{4.24}) and combining with Lemma~\ref{l2.3}, we obtain the result.~\qed \\
{\bf Remark:} By (\ref{3.9}) and (\ref{4.22}), we guess  $\nu_{n}^{\ast}=O(\frac{\log n}{n^{2}})$ when $\sigma =2$. However, the proof has not be obtained yet.

\noindent{\bf Proof of Theorem 3.4.}
For all $s> 0$, we have
\beqnn
P(Z_{n}\leq k_{n})
&=& P(e^{-\frac{s}{n}Z_{n}}\geq e^{-\frac{s}{n}k_{n}})\\
&\leq & E(e^{-\frac{s}{n}Z_{n}})e^{\frac{s}{n}k_{n}}\\
&=& H_{n}(e^{-\frac{s}{n}})e^{\frac{s}{n}k_{n}}.
\eeqnn
Letting $s=\frac{n}{k_{n}}$ and applying Proposition~\ref{p4.1}, we have
\beqnn
P(Z_{n}\leq k_{n})\leq C_{3}(1+\gamma \frac{n}{k_{n}})^{-\sigma}.
\eeqnn ~\qed

\noindent{\bf Proof of Theorem 3.5.}
Let $0< y_0 < R-1$, the sequence $y_{n}$ be defined by the equation
$$A(1+y_{n+1})=1+y_{n}.$$ It is not difficult to see that $A(1+y)\geq 1+y$ for $y\geq 0$. Therefore, the sequence $y_{n}$ decreases. Then,
\beqnn
\log H_{n}(1+y_{n})
&=&\sum\limits_{m=0}^{n-1}\log B(1+y_{n-m})\\
&=&\sum\limits_{i=1}^{n}\log B(1+y_{i})\\
&=&\sum\limits_{i=1}^{n}[B(1+y_{i})-1-\frac{1}{2\theta_{5}^{2}}(B(1+y_{i})-1)^{2}]\\
&=&\sum\limits_{i=1}^{n}[B^{'}(\theta_{6})y_{i}-\frac{B^{'}(\theta_{6})^{2}}{2\theta_{5}^{2}}y_{i}^{2}]\\
&\leq & B^{'}(1+y_{0})\sum\limits_{i=1}^{n}y_{i}+C_{11}\sum\limits_{i=1}^{n}y_{i}^{2},
\eeqnn
where $1\leq \theta_{5}\leq B(1+y_{i})$, $1\leq \theta_{6}\leq 1+y_{i},~i\geq 1$, and $C_{11}=-\frac{1}{2}(\frac{B'(1-)}{B(1+y_{0})})^{2}$. Therefore,
\beqnn
H_{n}(1+y_{n})\leq \exp\{B^{'}(1+y_{0})\sum\limits_{i=1}^{n}y_{i}+C_{11}\sum\limits_{i=1}^{n}y_{i}^{2}\}.
\eeqnn
By    \cite[Lemma 1]{N03},
\beqnn
\sum\limits_{i=0}^{n-1}y_{i}=\frac{1}{\gamma}\ln (1+\gamma ny_{0})+O(y_{0}),~~\sum\limits_{i=0}^{n-1}y_{i}^{2}=O(y_{0}).
\eeqnn
Setting $y_{0}=\frac{k_n}{\gamma^{2}n^{2}}-\frac{1}{\gamma n}$, we obtain
\beqlb\label{111}
H_{n}(1+y_{n})\leq (\frac{k_n}{\gamma n})^{B^{'}(1+\frac{k_n}{\gamma^{2}n^{2}}-\frac{1}{\gamma n})\frac{1}{\gamma}}.
\eeqlb
Note that for all $y> 0$,
\beqlb\label{222}
H_{n}(1+y)
&=& \sum\limits_{j=0}^{\infty}P(Z_{n}=j)(1+y)^{j}\nonumber\\
&\geq & \sum\limits_{j=k_n}^{\infty}P(Z_{n}=j)(1+y)^{j}\nonumber\\
&\geq & (1+y)^{k_n}P\{Z_{n}\geq k_n\}.
\eeqlb
According to \cite{N03},
\beqlb\label{333}
(1+y_{n})^{-k_n}=\exp\{-\frac{k_n}{\gamma n}+1-\frac{1}{\gamma}\lambda \frac{k_n}{n^{2}}\ln \frac{k_n}{n}\}(1+O(\frac{k_n}{n^{2}})).
\eeqlb
The theorem is proved  by combining (\ref{111})--(\ref{333}).
~\qed

\noindent{\bf Proof of Corollary  3.6.}
For every $t> 0$ we define $D_{n}(t)=e^{tZ_{n}},~n\geq 1$. It is easy to check that $\{D_n\}$ is a submartingale with respect to the natural $\sigma$-algebra generated by $\{Z_n\}$. By the Doob's inequality,
\beqnn
P(M_{n}\geq k)=P(\max\limits_{i\leq n}D_{i}(t)\geq e^{tk}) \leq \frac{ED_{n}(t)}{e^{tk}}=\frac{H_{n}(e^{t})}{e^{tk}}.
\eeqnn
Define $y_{n}$ as in the proof of Theorem~\ref{t3.5} and let $t=\log (1+y_{n})$. Then by the proof of Theorem~\ref{t3.5}, we can obtain the desired result.~\qed

~$
\end{math}

\end{document}